\documentclass[11pt,a4paper]{amsart} %mathkap} %
\usepackage[all]{xy}
\usepackage{graphics}
\usepackage{epigraph}
\usepackage{soul}

\usepackage{bbm}

\newcommand {\supplus}{\mathop{{\supset}\llap{\raise
0.5pt\hbox{\normalfont\small+}\hskip 0.5pt}}}
\newcommand {\subplus}{\mathop{{\subset}\llap{\raise
0.5pt\hbox{\normalfont\small+}\hskip 0.5pt}}}
\newcommand {\divby}  {\lower 0.15ex \hbox{\,\vdots\,}}
\newcommand {\Cee}    {{\mathbb  C}}

\newcommand {\Kee}    {{\mathbb  K}}

\newcommand {\Zee}    {{\mathbb  Z}}

\newcommand {\fag}    {{\mathfrak{ag}}}

\newcommand {\fder}   {{\mathfrak{der}}}   %

\newcommand {\fg}     {{\mathfrak{g}}}    %
\newcommand {\fh}     {{\mathfrak{h}}}

\newcommand {\fk}     {{\mathfrak{k}}}

\newcommand {\fm}     {{\mathfrak{m}}}

\newcommand {\fosp}   {{\mathfrak{osp}}}

\newcommand {\fpsl}   {{\mathfrak{psl}}}

\newcommand {\fs}     {{\mathfrak{s}}}

\newcommand {\fsi}    {{\mathfrak{si}}}
\newcommand {\fsl}    {{\mathfrak{sl}}}

\newcommand {\fsp}    {{\mathfrak{sp}}}

\newcommand {\fv}     {{\mathfrak{v}}}     %
\newcommand {\fvect}  {{\mathfrak{vect}}}   %

\newcommand {\fz}     {{\mathfrak{z}}}

\newcommand {\cal} {\mathcal}

\newcommand {\cD}     {{\cal D}}

\newcommand {\cO}     {{\cal O}}

\def \opname#1#2%
  {\expandafter\newcommand \csname #1\endcsname {{\mathop{#2}\nolimits}}}

\newcommand{\rmname}[1]
  {\expandafter\newcommand \csname #1\endcsname {{\operatorname{#1}}}}

\newcommand{\rmnameii}[2]
  {\expandafter\newcommand \csname #1\endcsname {{\operatorname{#2}}}}

\rmname{act} \rmname{Ad} \rmname{Add} \rmname{ad} \rmname{Alt}
\rmname{alt} \rmname{Ann} \rmname{antidiag} \rmname{Ber}
\rmname{ber} \rmname{Bil} \rmname{Br} \rmname{card} \rmname{ch}
\rmname{Char} \rmname{cem} \rmname{cj} \rmname{Cliff}
\rmname{cntr} \rmname{codim} \rmname{Coind} \rmname{const}
\rmname{col} \rmname{cork} \rmname{cpr} \rmname{diag}
 \rmnameii{Div}{div} \rmname{Def} \rmname{Deg}
\rmname{Der} \rmname{Diff} \rmname{Dim} \rmname{End} \rmname{Even}
\rmname{Ext} \rmname{gr} \rmname{Hom} \rmname{HT}
\rmnameii{Ht}{ht} \rmname{hwt} \rmname{Id} \rmname{id}
\rmname{ind} \rmname{Ind} \rmname{Inf} \rmname{irr} \rmname{Le}
\rmname{Lie} \rmname{lwt} \rmname{mult} \rmname{Mat} \rmname{Mor}
\rmname{nm} \rmname{Ob} \rmname{Odd} \rmname{Osc} \rmname{per}
\rmname{Pic} \rmname{pr} \rmname{pro} \rmname{Prime} \rmname{Proj}
\rmname{prt} \rmname{pt} \rmname{Q} \rmname{qet} \rmname{qtr}
\rmname{rd} \rmname{rk} \rmname{row} \rmname{Res} \rmname{salt}
\rmname{Sch} \rmname{SBr} \rmname{sdim}\rmname{scalar}
\rmname{Ser} \rmname{sign} \rmname{Smbl} \rmname{spin}
\rmname{ssym} \rmname{str} %\rmname{st}
\rmname{sgn} \rmname{sq}
\rmname{symm} \rmname{supp} \rmname{Supp} \rmname{St}
\rmname{Spec} \rmname{Spm} \rmname{tr} \rmname{vpt} \rmname{Vect}
\rmname{weyl} \rmname{Weyl} \rmname{Witt}

\opname{vvol}  {{v\hspace{-0.1ex}o\hspace{-0.02ex}l\/}}
\opname{pnt}  {\text{\normalfont pt}} \opname{Span} {{Span}}
\opname{slim} {\overline{\lim}} \opname{Vol}
{{V\hspace{-0.55ex}o\hspace{-0.02ex}l\/}} \opname{Par}
{{P\hspace{-0.3ex}a\hspace{-0.05ex}r\/}}

\newcommand {\ev} {{\bar0}}

\newcommand {\bcdot}   {\mathbin{\hbox{\raise.4ex\hbox{\bf.}}}} % bold \cdot

\newcommand {\secno} {}
\newcommand {\ssecfont} {\normalfont\bf}

\newtheorem*{Theorem}{\secno Theorem}

\newenvironment {th*}[1]
    {\gdef\thname{#1} \begin{thn}}%
    {\end{thn}}
\newtheorem*{thn} {\thname}

\theoremstyle{definition}
\newtheorem*{Example}{\secno Example}

\newenvironment {ex*}[1]
    {\gdef\thname{#1} \begin{exn}}%
    {\end{exn}}
\newtheorem*{exn}{\thname}

\theoremstyle{remark}

\newenvironment {rem*}[1]
    {\gdef\thname{#1} \begin{remn}}%
    {\end{remn}}
\newtheorem*{remn}{\thname}

\newcommand {\ssec}{\subsection*}

\newcommand {\ssbegin}[2]
  {\def \secno {\gdef \secno {}{\ssecfont #1. }}%
   \begin{#2}}

\begin{document}

\title[Simple Lie superalgebras]{Simple Lie superalgebras and
nonintegrable distributions in characteristic $p$}

\author{Sofiane Bouarroudj${}^1$, Dimitry Leites${}^2$}

\address{${}^1$Department of Mathematics, United Arab Emirates University, Al
Ain,
Po. Box: 17551; Bouarroudj.sofiane@uaeu.ac.ae\\
${}^2$MPIMiS, Inselstr. 22, DE-04103 Leipzig, Germany\\
on leave from Department of Mathematics, University of Stockholm,
Roslagsv. 101, Kr\"aft\-riket hus 6, SE-104 05 Stockholm,
Sweden; mleites@math.su.se, leites@mis.mpg.de}



\begin{abstract} Recently, Grozman and Leites returned to the original Cartan's
description of Lie algebras to interpret the Melikyan algebras
(for $p\leq 5$) and several other little-known simple Lie algebras
over algebraically closed fields for $p=3$ as subalgebras of Lie
algebras of vector fields preserving nonintegrable distributions
analogous to (or identical with) those preserved by G(2), O(7),
Sp(4) and Sp(10). The description was performed in terms of
Cartan-Tanaka-Shchepochkina prolongs using Shchepochkina's
algorithm and with the help of {\bf SuperLie} package. Grozman and
Leites also found two new series of simple Lie algebras.

Here we apply the same method to distributions preserved by one of
the two exceptional simple finite dimensional Lie superalgebras
over $\Cee$; for $p=3$, we obtain a series of new simple Lie
superalgebras and an exceptional one.
\end{abstract}

\thanks{We are thankful to P.~Grozman and I.~Shchepochkina for help;
DL is thankful to MPIMiS, Leipzig, for financial support and most
creative environment.}

%\date{Received May 1, 2006}

\maketitle

\epigraph{In memory of Felix Aleksandrovich Berezin}

\section{Introduction}

F. A. Berezin and supersymmetries are usually associated with
physics. However, Lie superalgebras --- infinitesimal
supersymmetries --- appeared in topology at approximately the same
time as the word \lq\lq spin" appeared in physics and it were
these examples that Berezin first had in mind. The Lie
superalgebras of topologists were often the ones over fields of
characteristic $p>0$, moreover, over finite fields, see \cite{W,
CL}. Since the natural symmetries are usually related to {\bf
simple} Lie (super)algebras, and the latter are easier to study
over algebraically closed fields, the attention of mathematicians
became focused on these, and here we follow the trend, although it
is clear that {\bf simple} Lie algebras and {\bf simple} Lie
superalgebras (for any $p$) are not as interesting as their \lq\lq
relatives" (central extensions, algebras of derivations, etc.).
Berezin was also interested in these topics but having felt that
H.~Galois's passionate words \lq\lq je n'ai pas de temps" were
applicable to him, he choose to concentrate on other things.
Lately, Lie algebras and superalgebras for $p$ prime return to the
scene as prime characters, see \cite{BKK, BKR}.

The habitual nowadays description of simple finite dimensional Lie
algebras in terms of abstract root systems, although enables one
to advance rather far in the study of representations, reminds us
V.~Arnold's description of Leibniz's contribution to Calculus:
\lq\lq Leibniz's notations are so good that using them everybody
can  nowadays study, and even teach, Calculus without any
understanding of what one is doing". In \cite{GL4}, Grozman and
Leites, keeping roots handy, applied Cartan's initial (now
practically forgotten) description of Lie algebras, not
necessarily simple ones, in terms of nonintegrable distributions
these algebras preserve. They interpreted a number of
ill-described simple Lie superalgebras over algebraically closed
fields of characteristic $p=5$ and 3, and discovered two (or
rather three) new series of simple Lie algebras for $p=3$. These
Lie algebras seemingly have no counterparts over $\Cee$ but
actually they do: they preserve the same nonintegrable
distributions as some of the classical simple Lie algebras (for
the lack of space, this interpretation was implicit in \cite{GL4};
the explicit description of the distributions in terms of Pfaff
equations can be easily obtained by means of Shchepochkina's
algorithm \cite{Shch} and {\bf SuperLie} package \cite{Gr}).

\ssec{1.1. Cartan's description of Lie algebras} The method of
constructing new Lie algebras used in \cite{GL4} is as follows
(all nonstandard terms will be described in the main text):

\medskip

(GL1) Take a simple finite dimensional Lie algebra $\fg$ over
$\Cee$, take its form over $\Zee$ and reduce all structure
constants modulo $p$. We may assume now that $\fg$ is being
considered over the ground field $\Kee$ of characteristic $p$.

(GL2) For each of the \lq\lq simplest" $\Zee$-gradings $\fg=\oplus
\fg_i$ of $\fg$, consider a version (adjusted to the algebra of
divided powers) of the Cartan-Tanaka-Shchepochkina (CTS) prolong
of the \lq\lq beginning" (see \cite{Shch}) part
$\fg=\mathop{\oplus}\limits_{i\leq
k_0}
\fg_i$, where $k_0$ is usually equal to 0 or 1.

(GL3) If the CTS prolong does not coincide with $\fg$, it contains
a simple ideal. Single out this ideal.

\medskip
The situation with the classification of simple finite dimensional
Lie algebras and Lie superalgebras is described in the next
subsection.

\ssec{1.2. The Kostrikin-Shafarevich conjecture and its super
counterpart} In \cite{KS}, Kostrikin and Shafarevich conjectured a
description of all simple restricted finite dimensional Lie
algebras $\fg$ over an algebraically closed field $\Kee$ of
characteristic $p>7$. Dropping the restrictedness condition, a
generalized KSh-conjecture states that all simple finite
dimensional Lie algebras can be obtained as the end product of the
following steps:

\medskip

(KS1) Take a $\Zee$-form $\fg_\Zee$ of a simple finite dimensional
Lie algebra $\fg$ over $\Cee$ and tensor this $\Zee$-form by
$\Kee$ over $\Zee$ and take all simple subquotients $\fsi(\fg)$ of
$\fg_\Zee\otimes_\Zee\Kee$ for all such $\fg$;

(KS2) take an analog $\fg_\Kee$ (with an algebra of divided powers
over $\Kee$ instead of polynomials) of a simple vectorial Lie
algebras $\fg$ over $\Cee$ and take all simple subquotients
$\fsi(\fg_\Kee)$ of $\fg_\Kee$ for all such $\fg$;

(KS3) take nontrivial deformations of the simple Lie algebras obtained at steps
(KS1)--(KS2), if any exists.

\medskip

This conjecture is now proven; it is even true for $p=7$. For
$p=5$, the list of simple finite dimensional Lie algebras contains
one more item: Melikyan\footnote{For their interpretation, see
\cite{GL4}.} algebras (and nothing else, as proven by Premet and
Strade, see \cite{S}). For $p=3$, there is no even feeling that a
complete list of simple algebras is obtained; for $p=2$, the mood
is even gloomier, cf. \cite{S}.

{\bf Super case}. Even when \cite{ALS} was being written (and the
classification of simple vectorial Lie superalgebras was not yet
even announced, see \cite{LSh, K3}), Leites and Shchepochkina
conjectured\footnote{Delivered at various seminars and conferences
starting from 1977, it was never published until now.} that {\bf
the same steps (KS1)--(KS3) applied to either simple finite
dimensional Lie superalgebras or the simple vectorial Lie
superalgebras (everything over $\Cee$) yield all simple finite
dimensional Lie superalgebras over an algebraically closed field
$\Kee$ of characteristic $p>7$, perhaps, even for $p=7$}. To prove
this conjecture is a challenging problem.

\ssec{1.3. Our result} Here we apply the method (GL1)--(GL3) to a
simple Lie superalgebra $\fag(2)$ and obtain new simple finite
dimensional Lie superalgebras for $p>2$. The novelty of our
examples is due to the fact that even the conjectural list of
simple Lie superalgebras (sec.~1.2) only covers $p>7$. The above
procedure (GL1)--(GL3) (description of deforms of the results thus
obtained should follow) is our method to get new examples of
simple Lie (super)algebras; to speak about completeness of the
bestiarium thus obtained is too early.

We start our quest for new simple Lie superalgebras with an
algebra related to $\fg(2)$\footnote{We denote the exceptional Lie
algebras in the same way as the serial ones, like $\fsl(n)$; we
thus avoid confusing $\fg(2)$ with the second component $\fg_2$ of
a $\Zee$-graded Lie algebra $\fg$.}, the latter being of interest
lately for physical applications \cite{AW}.

\section{Background}

For the background on Linear
Algebra in Superspaces and the list of simple Lie superalgebras, see \cite{LSh},
\cite{K2},
and also \cite{K3} and references therein.

\ssec{2.1. Integer bases in Lie superalgebras} Let $A=(a_{ij})$ be
an $n\times n$ matrix. A {\it Lie superalgebra $\fg=\fg(A)$ with
Cartan matrix} $A=(a_{ij})$, is given by its {\it Chevalley
generators}, i.e., elements $X^{\pm}_i$ of degree $\pm 1$ and $H_i
= [X_{i}^+, X_{i}^-]$ (of degree 0) that satisfy the relations
(hereafter in similar occasions either all superscripts $\pm$ are
$+$ or all are $-$)
\begin{equation}
\label{g(a)} {}[X_{i}^+, X_{j}^-]  = \delta_{ij}H_i, \quad [H_{i},
H_{j}] = 0, \quad [H_i, X_j^{\pm}] = \pm a_{ij}X_j^{\pm},
\end{equation}
and additional relations $R_i=0$ whose left sides are implicitly
described, for a general Cartan matrix, as
\begin{equation}
\label{myst}
\renewcommand{\arraystretch}{1.4}
\begin{array}{l}
\text{\lq\lq the $R_i$ that generate the maximal ideal $I$ such
that}\\
I\cap\Span(H_i \mid 1\leq i\leq n)=0. "
\end{array}
\end{equation}
For simple (finite dimensional) Lie algebras over $\Cee$, instead
of implicit description (\ref{myst}) we have an explicit
description (Serre relations) in terms of  Cartan matrices. For
definition of Cartan matrices of Lie superalgebras, see
\cite{GL1}. In particular, for all simple Lie superalgebras of the
form $\fg=\fg(A)$, except for the deforms of $\fosp(4|2)$, there exist
bases with respect to which all structure constants are integer.
In particular, this applies to $\fag(2)$, first discovered by
Kaplansky, see \cite{K1C, K2}. For presentations of $\fag(2)$, see
\cite{GL1}.

{\bf For vectorial Lie superalgebras},  integer bases are
associated with $\Zee$-forms of $\Cee[x]$ --- a supercommutative
superalgebra in $a$ (ordered for convenience) indeterminates $x =
(x_1,...,x_a)$ of which the first $m$ indeterminates are even and
the rest $n$ ones are odd ($m+n=a$). For a multi-index
$\underline{r}=(r_1, \ldots , r_a)$, we set
$$
u_i^{r_{i}} := \frac{x_i^{r_{i}}}{r_i!}\quad \text{and}\quad
u^{\underline{r}} := \prod\limits_{1\leq i\leq a} u_i^{r_{i}}.
$$
The idea is to formally replace fractions with $r_i!$ in
denominators by inseparable symbols $u_i^{r_{i}}$ which are
well-defined over fields of prime characteristic. Clearly,
\begin{equation}
\label{divp} u^{\underline{r}} \cdot u^{\underline{s}} = \binom
{\underline{r} + \underline{s}} {\underline{r}} u^{\underline{r} +
\underline{s}}, \quad \text{where}\quad\binom {\underline{r} +
\underline{s}} {\underline{r}}:=\prod\limits_{1\leq i\leq a}\binom
{r_{i} + s_{i}} {r_{i}}.
\end{equation}
For a set of positive integers $\underline{N} = (N_1,..., N_m)$,
denote
\begin{equation}
\label{u;N}
\renewcommand{\arraystretch}{1.4}
\begin{array}{l}
\cO(m; \underline{N}):=\Kee[u;
\underline{N}]:=\\
\Span_{\Kee}(u^{\underline{r}}\mid r_i < p^{N_{i}}\;\text{ for
$i\leq m$ and $r_i=0$ or 1 for $i>m$}).\end{array}
\end{equation}
As is clear from (\ref{divp}), $\Kee[u; \underline{N}]$ is a
subalgebra of $\Kee[u]$. The algebra $\Kee[u]$ and its subalgebras
$\Kee[u; \underline{N}]$ are called the {\it algebras of divided
powers}; they are analogs of the algebra of functions.

An important feature that differs $\Kee[u]$  from $\Kee[u;
\underline{N}]$ is the number of generators: whereas $\Kee[u]$ is
generated by the $u_i$, the algebra $\Kee[u; \underline{N}]$ has,
additionally, the generators $u_i^{(p^{k_i})}$ for every
$k_i$ such that
$1<k_i<N_i$. Since any derivation of a given algebra is completely
determined by its value on every generator of the algebra, the Lie
algebra of all derivations of $\Kee[u; \underline{N}]$ is much
larger than the Lie algebra of {\it special derivations} whose
generators behave like partial derivatives:
\begin{equation}
\label{special}
\partial_i(u_j^{(k)})=\delta_{ij}u_j^{(k-1)}.
\end{equation}
In what follows, we only consider special derivations, e.g., in
(\ref{specialW}).

The simple vectorial Lie algebras over $\Cee$ have only one
parameter: the number of indeterminates. If $\Char \Kee =p>0$, the
vectorial Lie algebras acquire one more parameter:
$\underline{N}$. For Lie superalgebras, $\underline{N}$ only
concerns the even indeterminates. Let
\begin{equation}
\label{specialW} \fvect (m; \underline{N}|n) \;\text{ a.k.a }W(m;
\underline{N}|n)\; :=\fder \Kee[u; \underline{N}]
\end{equation}
be the general vectorial Lie algebra.

\ssec{2.2. $\Zee$-gradings} Recall that every $\Zee$-grading of a
given vectorial algebra is determined by setting $\deg u_i=r_i\in
\Zee$; every $\Zee$-grading of a given Lie superalgebra $\fg(A)$
is determined by setting $\deg X^{\pm}_i=\pm r_i\in \Zee$.

For the Lie algebras of the form $\fg(A)$, we set
\begin{equation}\label{selec}
\deg X^{\pm}_i=\pm\delta_{i, i_j}\;\text{ for any $i_j$ from a
selected set $\{i_1, \dots , i_k\}$}
\end{equation}
and say that we have \lq\lq {\sl selected}\rq\rq certain $k$
Chevalley generators (or respective nodes of the Dynkin graph).
Yamaguchi's theorem cited below shows that, in the study of Cartan
prolongs defined below, the first gradings to consider are the
ones with all $r_i=0$ except $k$ of them ($1\leq k\leq 2$) \lq\lq
{\sl selected}\rq\rq Chevalley generators for which $r_i=1$. In
this paper we consider the simplest gradings, for $k=1$.

For vectorial algebras, filtrations are more natural than
gradings; the very term \lq\lq vectorial" means, actually, that
the algebra is endowed with  a particular ({\it Weisfeiler})
filtration, see \cite{LSh}. Unlike Lie algebras, the vectorial Lie
{\sl super}algebras can sometimes be regraded into each other;
various realizations as vectorial algebras are described by means
of one more parameter --- regrading $\underline{r}$ --- with a
\lq\lq {\it standard grading}" as a point of reference:
\begin{equation}\label{regr}
\renewcommand{\arraystretch}{1.4}
\begin{array}{l}
\fvect (m; \underline{N}|n; \underline{r}) \;\text{ a.k.a }W(m;
\underline{N}|n; \underline{r})\; :=\fder \Kee[u; \underline{N}],
\text{ where }\\
\deg u_i = r_i\;\text{  is a grading of $\cO(m; \underline{N}|n)$}.\end{array}
\end{equation}
For $W(m; \underline{N}|n)$, the standard grading is
$\underline{r}=(1, \dots , 1)$. For the contact algebras
$\fk(2n+1, \underline{N})$ that preserve the Pfaff equation
$\alpha(X)=0$ for $X\in\fvect(2n+1|m)$, where (see \cite{Le})
\begin{equation}\label{cont}
\alpha=\begin{cases}
dt-\mathop{\sum}\limits_{i\leq n}(p_idq_i-q_idp_i)+
\sum\limits_{j\leq m}\theta_jd\theta_j&\text{if $p\neq 2$,}\\
dt+\sum\limits_{1\leq i\leq k}
x_idx_{k+i}\begin{cases}&\text{for $n=2k$ and $x_1,\dots,x_n$ all
even or all odd}\\
+x_ndx_n&\text{for $n=2k+1$ and $x_1,\dots x_n$ odd}\end{cases}
&\text{if $p= 2$}.\end{cases}
\end{equation}
the standard grading is $\deg t=2$ and the degree all other
indeterminates being equal to $1$.

\ssec{2.3. Cartan prolongs} Let $\fg_0$ be a Lie algebra,
$\fg_{-1}$ a $\fg_0$-module. Let us define  the $\Zee$-graded Lie
algebra $(\fg_{-1}, \fg_{0})_*=\oplus_{i\geq -1}\fg_i$  called the
{\it complete Cartan prolong} (the result of the Cartan {\it
prolongation}) of the pair $(\fg_{-1}, \fg_{0})$. Geometrically
the Cartan prolong is the maximal Lie algebra of symmetries of the
$G$-structure (here: $\fg_0=\text{Lie}(G)$) on $\fg_{-1}$. The
components $\fg_i$ for $i>0$ are defined recursively.

First, recall that, for any (finite dimensional)
vector space $V$, we have
$$
\Hom(V, \Hom(V,\ldots, \Hom(V,V)\ldots)) \simeq L^{i}(V, V,
\ldots, V; V),
$$
where $L^{i}$ is the space of $i$-linear maps and we have
$(i+1)$-many $V$'s on both sides. Now, we recursively define, for
any $v_1, \dots, v_{i+1}\in \fg_{-1}$ and any $i > 0$:
$$
\renewcommand{\arraystretch}{1.4}
\begin{array}{ll}
\fg_i =& \{X\in \Hom(\fg_{-1}, \fg_{i-1})\mid X(v_1)(v_2, v_3,
...,
v_{i+1}) =\\
&X(v_2)(v_1, v_3, ..., v_{i+1})\}.
\end{array}
$$

Let the $\fg_0$-module $\fg_{-1}$ be faithful. Then, clearly,
$$
(\fg_{-1}, \fg_{0})_*:=\oplus\fg_i\subset \fvect (m) = \fder~
\Kee[x_1, \ldots , x_m], \; \text{ where}\; m = \dim~ \fg_{-1}.
$$
Moreover, setting $\deg x_i=1$ for all $i$, we see that
$$
\renewcommand{\arraystretch}{1.4}
\begin{array}{l}
\fg_i = \{X\in \fvect (m)\mid \deg X=i,\; [X, \partial]\in
\fg_{i-1}\;\text{ for any }\partial\in \fg_{-1}\}.
\end{array}
$$

Now it is subject to an easy verification that the Cartan prolong
$(\fg_{-1}, \fg_{0})_*$ forms a subalgebra of $\fvect (m)$. (It is
also easy to see that $(\fg_{-1}, \fg_{0})_*$  is a Lie algebra
even if $\fg_{-1}$ is not a faithful $\fg_0$-module; but then it
can not be realized as a subalgebra of $\fvect(m)$.)

Obviously, for $p>0$, there is a series of Cartan prolongs
labelled by $\underline{N}$ and the same applies to the
CTS-prolongs and the partial prolongs defined in the next
subsections.

\ssec{2.4. Nonholonomic manifolds. Cartan-Tanaka-Shchepochkina
(CTS) prolongs} Let $M^n$ be an $n$-di\-men\-si\-o\-nal manifold.
Recall that a {\it distribution} $\cD$ on $M$ is any subbundle of
the tangent bundle; $\cD$ is said to be {\it integrable} (in a
neighborhood $U$ of a point $m\in M$) if, for each point $x\in U$,
there is a local submanifold $S$ of $U$ containing $x$ and such
that the tangent bundle to $S$ is equal to $\cD$ restricted to
$S$. A criterion due to Frobenius states that $\cD$ is integrable
if and only if its sections form a Lie algebra. Let  $\cD$ be a
nonintegrable distribution; then there exists a sequence of strict
inclusions
$$
\cD= \cD_{-1}\subset \cD_{-2} \subset \cD_{-3} \dots \subset
\cD_{-d},
$$
where the fiber of $\cD_{-i}$ at a point $x\in M$ is
$$
\cD_{-i+1}(x)+ [\cD_{-1}, \cD_{-i+1}](x)
$$
(here $[\cD_{-1}, \cD_{-i-1}]=\Span\left([X, Y]\mid X\in \Gamma(\cD_{-1}), Y\in
\Gamma(\cD_{-i-1})\right)$) and $d$ is the least number such that
$$
\cD_{-d}(x)+[\cD_{-1}, \cD_{-d}](x) = \cD_{-d}(x).
$$
In case $\cD_{-d} = TM$ the distribution is called {\it completely
nonholonomic}. The number $d = d(M)$ is called the {\it
nonholonomicity degree}. A manifold $M$ with a distribution $\cD$ on
it will be referred to as {\it nonholonomic} one if $d(M)\neq 1$.
Let
\begin{equation}\label{n_i}
n_i(x) = \dim \cD_{-i}(x); \qquad n_0(x)=0; \qquad
n_d(x)=n-n_{d-1}.
\end{equation}
The distribution $\cD$ is said to be
{\it regular} if all the dimensions $n_i$ are constants on $M$. We
will only consider regular, completely nonholonomic distributions, and,
moreover,
satisfying certain transitivity condition (\ref{tr}) introduced below.

To the tangent bundle over a nonholonomic manifold $(M, \cD)$ we
assign a bundle of $\Zee$-graded nilpotent Lie algebras as
follows. Fix a point $pt\in M$. The usual adic filtration by
powers of the maximal ideal $\fm:=\fm_{pt}$ consisting of
functions that vanish at $pt$ should be modified because distinct
coordinates may have distinct \lq\lq degrees". The distribution
$\cD$ induces the following filtration in $\fm$:
\begin{equation}\label{2.1}
\renewcommand{\arraystretch}{1.4}
\begin{array}{ll}
\fm_k=&\{f\in\fm\mid X_1^{a_1}\ldots X_n^{a_n}(f)=0\;\text{ for
any $X_1, \dots,X_{n_1}\in \Gamma(\cD_{-1})$,  }\\
&  \text{$X_{n_1+1}, \dots, X_{n_2}\in \Gamma(\cD_{-2})$,\dots,
$X_{n_{d-1}+1}, \dots, X_{n}\in \Gamma(\cD_{-d})$}\\
& \text{ such that }\; \mathop{\sum}\limits_{1\leq i\leq d}\quad
\left(i\mathop{\sum}\limits_{n_{i-1}< j\leq n_{i}} a_j\right )\leq k\},
\end{array}
\end{equation}
where $\Gamma(\cD_{-j})$ is the space of germs at $pt$ of sections
of the bundle $\cD_{-j}$. Now, to a filtration
$$
\cD= \cD_{-1}\subset \cD_{-2} \subset \cD_{-3} \dots \subset \cD_{-d}=TM,
$$
we assign the associated graded bundle
$$
\gr(TM)=\oplus\gr\cD_{-i},\;\text{ where $\gr\cD_{-i}=\cD_{-i}/\cD_{-i+1}$}
$$
and the bracket of sections of $\gr(TM)$ is, by definition, the
one induced by bracketing vector fields, the sections of $TM$. We assume
a \lq\lq transitivity condition\rq\rq: The
Lie algebras
\begin{equation}\label{tr}
\gr(TM)|_{pt}
\end{equation}
induced at each point $pt\in M$ are isomorphic.

The grading of the coordinates $(\ref{2.1})$ determines a
nonstandard grading of $\fvect(n)$ (recall (\ref{n_i})):
\begin{equation}\label{gr}
\renewcommand{\arraystretch}{1.4}
\begin{array}{l}
\deg x_1=\ldots =\deg x_{n_1}=1,\\
\deg x_{n_1+1}=\ldots =\deg x_{n_2}=2,\\
\dotfill \\
\deg x_{n-n_{d-1}+1}=\ldots =\deg x_{n}=d.
\end{array}%\eqno{(2.9)}
\end{equation}
Denote by $\fv=\mathop{\oplus}\limits_{i\geq -d}\fv_i$ the algebra
$\fvect(n)$ with the grading $(\ref{gr})$. One can
show that the \lq\lq complete prolong'' of $\fg_-$ to be defined
shortly, i.e., $(\fg_-)_*:=(\fg_-, \tilde \fg_0)_*\subset \fv$, where $\tilde
\fg_0:=\fder_0\fg_-$, preserves $\cD$.

For nonholonomic manifolds, an analog of the group $G$ from the
term ``$G$-structure'', or rather of its Lie algebra,
$\fg=\text{Lie}(G)$, is the pair $(\fg_-, \fg_0)$, where $\fg_0$
is a  subalgebra of the $\Zee$-grading preserving Lie algebra of
derivations of $\fg_-$, i.e., $\fg_0 \subset \fder_0\,\fg_-$. If
$\fg_0$ is not explicitly indicated, we assume that $\fg_0
=\fder_0\,\fg_-$, i.e., is the largest possible.

Given a pair $(\fg_-, \fg_0)$ as above, define its {\it
Tanaka-Shchepochkina prolong} to be the maximal subalgebra
$(\fg_-, \fg_0)_*=\mathop{\oplus}\limits_{k\geq -d} \fg_k$ of
$\fv$ with given non-positive part $(\fg_-, \fg_0)$. For an
explicit construction of the components, see  \cite{Shch}.
If $\fg_-=\fg_{-1}$ the
Tanaka-Shchepochkina prolong turns into its particular case, the well-known
Cartan
prolong.

\ssec{2.5. Partial prolongs and projective structures} Let
$(\fg_-, \fg_0)_*$ be a depth $d$ Lie algebra; $\fh_1\subset
\fg_1$ be a $\fg_0$-submodule  such that $[\fg_{-1},
\fh_1]=\fg_0$. If such $\fh_1$ exists, define the $i$th partial
prolong of $(\mathop{\oplus}\limits_{i\leq 0}\fg_i, \fh_1)$ for
$i\geq 2$ to be
\begin{equation}
\label{partprol} \fh_{i}=\{D\in\fg_{i}\mid [D, \fg_{-1}]\in
\fh_{i-1}\}.
\end{equation}
Set $\fh_i=\fg_i$ for $i\leq 0$ and call
$\fh_*=\mathop{\oplus}\limits_{i\geq -d}\fh_i$ the Shchepochkina
{\it partial prolong} of $(\mathop{\oplus}\limits_{i\leq 0}\fg_i,
\fh_1)$, see \cite{Sh14, Shch}. (Of course, the partial prolong can also be
defined
if
$\fh_0$ is contained in $\fg_0$.)

\begin{Example} The $SL(n+1)$-action on the projective space
$P^n$ gives the embedding $\fsl(n+1)\subset\fvect(n)$; here
$\fsl(n+1)$ is a partial prolong of $\fvect(n)_{i\leq
0}\oplus\fh_1$ for some $\fh_1$.
\end{Example}

\ssec{2.6. Yamaguchi's theorem} Let
$\fs=\mathop{\oplus}\limits_{i\geq -d}\fs_i$ be a simple finite
dimensional Lie algebra. Let $(\fs_-)_*=(\fs_-, \fg_0)_*$ be the
CTS prolong with the maximal possible
$\fg_0=\fder_0(\fs_-)$.

\begin{Theorem}[\cite{Y}] Over $\Cee$, the isomorphism $(\fs_-)_*\simeq\fs$
holds
almost always. The exceptions (cases where
$\fs=\mathop{\oplus}\limits_{i\geq -d}\fs_i$ is a partial prolong
in $(\fs_-)_*=(\fs_-, \fg_0)_*$) are

{\em 1)} $\fs$ with the grading of depth $d=1$ (in which case
$(\fs_-)_*=\fvect(\fs_-^*)$);

{\em 2)} $\fs$ with the grading of depth $d=2$ and
$\dim\fs_{-2}=1$, i.e., with the \lq\lq contact'' grading, in
which case $(\fs_-)_*=\fk(\fs_-^*)$ (these cases correspond to
\lq\lq selection" of the nodes on the Dynkin graph connected with
the node for the maximal root on the extended graph);

{\em 3)} $\fs$ is either $\fsl(n+1)$ or $\fsp(2n)$ with the
grading determined by \lq\lq selecting'' the first and the $i$th
of simple coroots, where $1<i<n$ for $\fsl(n+1)$ and $i=n$ for
$\fsp(2n)$. (Observe that $d=2$ with $\dim \fs_{-2}>1$ for
$\fsl(n+1)$ and $d=3$ for $\fsp(2n)$.)

Moreover, the isomorphism  $(\fs_-, \fs_0)_*\simeq\fs$ also holds
almost always. The cases where it fails (the ones where a
projective action is possible) are $\fsl(n+1)$ or $\fsp(2n)$ with
the grading determined by \lq\lq selecting'' only one (the first)
simple coroot. \end{Theorem}

\section{Our Examples}

Our examples, as well as Melikyan's ones, and those of \cite{GL4},
are due to the fact that, for $p$ small, the analog of Yamaguchi's
theorem is false both for Lie algebras and Lie superalgebras.

Serganova and van de Leur showed (\cite{Se}, \cite{vdL};
\cite{Se1}) that the Lie superalgebra $\fg=\fag(2)$ has the
following four non-equivalent Cartan matrices:
$$
1)\; \begin{pmatrix} 0 & 1 & 0 \\ -1 & 2 & -3 \\ 0 & -1 & 2
\end{pmatrix}\quad 2)\; \begin{pmatrix}
0 & 1 & 0 \\ -1 & 0 & 3 \\ 0 & -1 & 2
\end{pmatrix}\quad 3)\; \begin{pmatrix}
0 & -3 & 1 \\ -3 & 0 & 2 \\ -1 & -2 & 2
\end{pmatrix}\quad 4) \; \begin{pmatrix}
2 & -1 & 0 \\ -3 & 0 & 2 \\ 0 & -1 & 1
\end{pmatrix}
$$
Here we only consider the simplest $\Zee$-grading $r$ and only the first matrix.

\ssec{4.1. The first Cartan matrix and $r=(1,0,0)$} Then
$\fg=\mathop{\oplus}\limits_{|i|\leq 2}\fg_i$, where
$\dim\fg_{-2}=1$ and $\sdim \fg_{-1}=0|7$. Therefore,
$\fg\subset\fk(1|7)$. Observe that $\fg_0=\fg(2)\oplus\fz$, where
$\fz=\Span(t)$ is the center of $\fg_0$. Hereafter the elements of
the Lie superalgebra of contact vector fields are given in terms
of their generating functions in indeterminates introduced in
(\ref{alpha}).

From an explicit description of $\fg(2)$ and its first fundamental
representation in \cite{FH} we deduce an explicit form of the
non-positive elements of $\fg$ which we give in terms of the
generating functions with respect to the contact bracket
corresponding to the contact form
\begin{equation}\label{alpha}
\alpha=dt-\sum\limits_{i=1,
3, 4}(v_idw_i+w_idv_i)+2udu,
\end{equation}
where the $v_i$, $w_i$ and $u$
are odd and notations match \cite{FH}, p. 354, while our $X_i^+$
and $X_i^-$ correspond to $X_i$ and $Y_i$ of \cite{FH}, p.~340, respectively.

We also set $X_3^\pm :=[X_1^\pm , X_2^\pm ]$, $X_4^\pm
:=[X_1^\pm , X_3^\pm ]$, $X_5^\pm :=[X_1^\pm , X_4^\pm ]$, $X_6^\pm :=[X_2^\pm , X_5^\pm ]$.

To describe the $\fg_0$-module $\fg_{-1}=\Span(u\;\text{ and $v_i,
w_i$ for $i=1, 3, 4$})$, only the highest weight vector suffices:
$$
\renewcommand{\arraystretch}{1.4}
\begin{tabular}{|l|l|}
\hline
$\fg_{i}$&the generating functions of generators, as $\fg_{0}$-modules \\
\hline \hline
$\fg_{-2}$&$1$\\
\hline
$\fg_{-1}$& $v_4$\\
\hline
$\fg_0$&$\fz=\Span(t)$\\
&$X_1^+= -v_4w_3-uv_1$\\
&$X_2^+=v_3w_1$\\
&$X_1^-=-v_3w_4-uw_1$\\
&$X_2^-=v_1w_3$\\
\hline
\end{tabular}
$$
As expected,  for $p=0$
and $p>3$, the CTS prolong is isomorphic to $\fag(2)$.

For $p=3$, the Lie algebra $\fg_0$ is not simple, but has a simple
Lie subalgebra isomorphic to $\fpsl(3)$ generated by
$x_1^\pm=X_1^\pm$, and $x_2^\pm =[X_{1}, X_2^\pm]$ and spanned by
$\{X^\pm_1, X^\pm_3, X^\pm_4, H_1 \}$. Let
$\widetilde\fg_0:=\fpsl(3)\oplus\fz$, where $\fz=\Span(t +
v_{1}w_{1} +  v_{3} w_{3} + 2\,  v_{4} w_{4})$ is the center of
$\widetilde\fg_0$.

The $\widetilde\fg_0$-module $\fg_{1}$ splits into two irreducible
components: A $0|1$-dimensional, and $0|7$-dimensional with lowest
weight vectors, respectively:
$$
\renewcommand{\arraystretch}{1.4}
\begin{tabular}{|l|l|}
\hline
$\fg_{1}$&generating function\\
\hline \hline
$V_{1}'$&$v_1v_3w_4+v_1uw_1+v_3uw_3+2v_4uw_4+v_4w_1w_3$\\
$V_{1}''$& $tw_4+v_1w_1w_4+v_3w_3w_4+uw_1w_3$\\
\hline
\end{tabular}
$$

Since $\fg_{1}$ generates the positive part of the CTS prolong,
$[\fg_{1}, \fg_{-1}]=\widetilde\fg_{0}$, the
$\widetilde\fg_0$-module $\fg_{-1}$ is irreducible, and
$[\fg_{-1}, \fg_{-1}]=\fg_{-2}$, the standard criterion for
simplicity (\cite{S}) ensures that the CTS prolong is simple. Since
none of the finite dimensional simple Lie superalgebras over
$\Cee$ has grading of this form, except vectorial ones, and none
of vectorial simple Lie superalgebras $\fg$ have the simple part
of $\widetilde\fg_0$ isomorphic to $\fpsl(3)$ (which only exists
for $p=3$), we conclude that this Lie superalgebra is a new simple
Lie superalgebra indigenous to $p=3$. We denote it
$\text{Bj}(1;N|7)$.

The positive components of $\text{Bj}(1;N|7)$ are all of dimension
8 (the direct sums $\fg_k=\fg_k'\oplus\fg_k''$ of irreducible
$\fg_0$-modules of dimension 1 and 7, up to parity), except the ones of
the highest degree ($2(3^N-1)+3-2=2\cdot 3^N-1$), which are all of dimension 1,
and the second highest degree, which are all of dimension 7. Let $V_{k}'$ and
$V_{k}''$ be the lowest weight vectors (with respect to $\fg_0$) of $\fg_k'$ and $\fg_k''$,
respectively:
$$
\renewcommand{\arraystretch}{1.4}
\begin{tabular}{|l|l|}
\hline
$\fg_{k}$&generating function: $N=1$\\
\hline \hline
$V_{2}'$&$t^2+2v_4uw_1w_3+v_3v_4w_3w_4+v_1v_4w_1w_4+2v_1v_3w_1w_3+2v_1v_3uw_4$\\
$V_{2}''$& $2v_1uw_1w_4+2v_3uw_3w_4+2v_4w_1w_3w_4+tuw_4+2tw_3w_4$\\
\hline
$V_{3}'$&$tv_1v_3w_4+tv_1uw_1+tv_3uw_3+2tv_4uw_4+tv_4w_1w_3$\\
$V_{3}''$&
$t^2w_4+tuw_1w_3+2v_1v_3w_1w_3w_4+tv_1w_3w_4+tv_3w_3w_4+2v_4uw_3w_4$\\
\hline
$V_{4}''$&
$2v_1v_3uw_1w_3w_4+2tv_1uw_1w_4+2tv_3uw_3w_4+2tv_4w_1w_3w_4+t^2uw_4+2t^2w_1w_3$\\
\hline
$V_{5}'$&$2v_1v_3v_4uw_1w_3w_4+2t^2v_1v_3w_4+2t^2v_1uw_1+
2t^2v_3uw_3+t^2v_4uw_4+2t^2v_4w_1w_3
$\\
\hline
\end{tabular}
$$
Similarly (the components up to 3 being the same as for $N=1$):
$$
\renewcommand{\arraystretch}{1.4}
\begin{tabular}{|l|l|}
\hline
$\fg_{k}$&generating function: $N=2$\\
\hline
\hline
$V_{4}'$&$t^3+2tv_1v_3uw_4+2tv_1v_3w_1w_3+
tv_1v_4w_1w_4+tv_3v_4w_3w_4+2tv_4uw_1w_3$\\
$V_{4}''$&
$2v_1v_3uw_1w_3w_4+2tv_1uw_1w_4+2tv_3uw_3w_4+
2tv_4w_1w_3w_4+t^2uw_4+2t^2w_1w_3$\\
\hline
$V_{5}'$&$2v_1v_3v_4uw_1w_3w_4+2t^2v_1v_3w_4+
2t^2v_1uw_1+2t^2v_3uw_3+t^2v_4uw_4+2t^2v_4w_1w_3$\\
$V_{5}''$&
$t^3w_4+t^2v_1w_1w_4+t^2v_3w_3w_4+t^2uw_1w_3+2tv_1v_3w_1w_3w_4+2tv_4uw_1w_3w_4$\\
\hline
$V_{6}'$&$t^3uw_4+2t^3w_1w_3+2t^2v_1uw_1w_4+
2t^2v_3uw_3w_4+2t^2v_4w_1w_3w_4+2tv_1v_3uw_1w_3w_4$\\
$V_{6}''$&
$t^4+2t^2v_1v_3uw_4+2t^2v_1v_3w_1w_3+t^2v_1v_4w_1w_4+
t^2v_3v_4w_3w_4+2t^2v_4uw_1w_3$\\
\hline
$\dots$&$\dotfill$\\
\hline
$V_{16}''$&
$t^8uw_4+2t^8w_1w_3+2t^7v_1uw_1w_4+
2t^7v_3uw_3w_4+2t^7v_4w_1w_3w_4+2t^6v_1v_3uw_1w_3w_4$\\
\hline
$V_{17}'$&$2t^8v_1v_3w_4+2t^8v_1uw_1+2t^8v_3uw_3+
t^8v_4uw_4+2t^8v_4w_1w_3+2t^6v_1v_3v_4uw_1w_3w_4$\\
\hline
\end{tabular}
$$
Similarly (the components up to 15 being the same as for $N=2$):
$$
\renewcommand{\arraystretch}{1.4}
\begin{tabular}{|l|l|}
\hline
$\fg_{k}$&generating function: $N=3$\\
\hline \hline
$V_{16}'$&$t^8uw_4+2t^8w_1w_3+2t^7v_1uw_1w_4+2t^7v_3uw_3w_4
+2t^7v_4w_1w_3w_4+2t^6v_1v_3uw_1w_3w_4$\\
$V_{16}''$& $t^9+2t^7v_1v_3uw_4+2t^7v_1v_3w_1w_3
+t^7v_1v_4w_1w_4+t^7v_3v_4w_3w_4+2t^7v_4uw_1w_3$\\
\hline
$V_{17}'$&$2t^8v_1v_3w_4+2t^8v_1uw_1+2t^8v_3uw_3+t^8v_4uw_4
+2t^8v_4w_1w_3+2t^6v_1v_3v_4uw_1w_3w_4$\\
$V_{17}''$& $t^9w_4+t^8v_1w_1w_4+t^8v_3w_3w_4+t^8uw_1w_3
+2t^7v_1v_3w_1w_3w_4+2t^7v_4uw_1w_3w_4$\\
\hline
$\dots$&$\dotfill$\\
\hline
$V_{52}''$& $t^{26}uw_4+2t^{26}w_1w_3+2t^{25}v_1uw_1w_4+2t^{25}v_3uw_3w_4+
2t^{25}v_4w_1w_3w_4+2t^{24}v_1v_3w_1w_3w_4$\\
\hline
$V_{53}'$&$t^{26}v_1v_3w_4+t^{26}v_1uw_1+t^{26}v_3uw_3+t^{26}v_4uw_4
+t^{26}v_4w_1w_3+t^{24}v_1v_3v_4uw_1w_3w_4$\\
\hline
\end{tabular}
$$
Now, let us consider subalgebras of $\text{Bj}(1;N|7)$:

(i) Let $\fg_1'$ be generated by $V_1'$, as
$\widetilde\fg_0$-module. The component of degree 2 of the partial
CTS prolong $(\fg_-, \fg_0, \fg_1')_*$ is $0$, so the corresponding Lie algebra
is
not
simple.

(ii) Let $\fg_1''$ be generated by  $V_1''$, as
$\widetilde\fg_0$-module. Let $\text{bj}=(\fg_-, \fg_0, \fg_1'')_*$ denote the
partial CTS
prolong , then
$$
\dim \fg_2=\dim \Span(t^2+2v_1v_3uw_3+2v_1v_3w_1w_3
+v_1v_3w_1w_4+v_3v_4w_3w_4+2v_4uw_1w_3)=1
$$
and
$\fg_3=0$. We see that $[\fg_1'', \fg_{-1}]=\widetilde\fg_0$.
The corresponding partial CTS prolong  is of dimension
$(10|14)$ and the criteria for simplicity imply that $\text{bj}$
is simple. Since none of the known Lie superalgebra has such
$\text{bj}_\ev\simeq \fsl(2)\oplus \fpsl(3)$,   this simple Lie superalgebra is
new. Since the partial CTS prolong for $N=2$ is the same, it is an
exceptional Lie superalgebra.

Let us summarize:

\ssbegin{4.2}{Theorem} For the first Cartan matrix of $\fag(2)$
and its $\Zee$-grading $r=(100)$, the CTS prolongs return
$\fag(2)$, except for the case $p=3$, which yield the simple
subalgebras of series {\em $\text{Bj}(1;N|7)$}, and the
exceptional algebra {\em $\text{bj}$}.
\end{Theorem}

\ssec{Post scriptum} We do not regularly check arXiv news, and
missed several interesting papers by Elduque cited in \cite{CE}
and almost missed (the deadline for Berezin's volume being April
15) \cite{CE} (which appeared on May 15), where several new simple
Lie superalgebras for $p=3$, also related to $\fg(2)$, are
described by an approach different from ours. Elduque's
superalgebras look non-isomorphic to ours.

\end{document}